\documentclass[12pt]{article}
\usepackage{latexsym, amssymb}
\usepackage{dsfont}

\DeclareMathAlphabet\gothic{U}{euf}{m}{n}

\makeatletter
\def\eqnarray{\stepcounter{equation}\let\@currentlabel=\theequation
\global\@eqnswtrue
\tabskip\@centering\let\\=\@eqncr
$$\halign to \displaywidth\bgroup\hfil\global\@eqcnt\z@
  $\displaystyle\tabskip\z@{##}$&\global\@eqcnt\@ne
  \hfil$\displaystyle{{}##{}}$\hfil
  &\global\@eqcnt\tw@ $\displaystyle{##}$\hfil
  \tabskip\@centering&\llap{##}\tabskip\z@\cr}

\def\endeqnarray{\@@eqncr\egroup
      \global\advance\c@equation\m@ne$$\global\@ignoretrue}

\def\@yeqncr{\@ifnextchar [{\@xeqncr}{\@xeqncr[5pt]}}
\makeatother

\begin{document}
\bibliographystyle{tom}

\newtheorem{lemma}{Lemma}[section]
\newtheorem{thm}[lemma]{Theorem}
\newtheorem{cor}[lemma]{Corollary}
\newtheorem{voorb}[lemma]{Example}
\newtheorem{rem}[lemma]{Remark}
\newtheorem{prop}[lemma]{Proposition}

\newenvironment{remarkn}{\begin{rem} \rm}{\end{rem}}
\newenvironment{exam}{\begin{voorb} \rm}{\end{voorb}}

\newcommand{\gota}{\gothic{a}}
\newcommand{\gotb}{\gothic{b}}
\newcommand{\gothh}{\gothic{h}}
\newcommand{\gotq}{\gothic{q}}
\newcommand{\gott}{\gothic{t}}

\newcounter{teller}
\renewcommand{\theteller}{(\alph{teller})}
\newenvironment{tabel}{\begin{list}%
{\rm  (\alph{teller})\hfill}{\usecounter{teller} \leftmargin=1.1cm
\labelwidth=1.1cm \labelsep=0cm \parsep=0cm}
                      }{\end{list}}

\newcounter{tellerr}
\renewcommand{\thetellerr}{(\roman{tellerr})}
\newenvironment{tabeleq}{\begin{list}%
{\rm  (\roman{tellerr})\hfill}{\usecounter{tellerr} \leftmargin=1.1cm
\labelwidth=1.1cm \labelsep=0cm \parsep=0cm}
                         }{\end{list}}

\newcommand{\Ni}{\mathds{N}}
\newcommand{\Qi}{\mathds{Q}}
\newcommand{\Ri}{\mathds{R}}
\newcommand{\Ci}{\mathds{C}}
\newcommand{\Ti}{\mathds{T}}
\newcommand{\Zi}{\mathds{Z}}
\newcommand{\Fi}{\mathds{F}}

\newcommand{\proof}{\mbox{\bf Proof} \hspace{5pt}} 
\newcommand{\vertspace}{\vskip10.0pt plus 4.0pt minus 6.0pt}

\newcommand{\RRe}{\mathop{\rm Re}}
\newcommand{\IIm}{\mathop{\rm Im}}

\newcommand{\spann}{\mathop{\rm span}}

\hyphenation{groups}
\hyphenation{unitary}

\newcommand{\tfrac}[2]{{\textstyle \frac{#1}{#2}}}

\newcommand{\cd}{{\cal D}}
\newcommand{\ch}{{\cal H}}
\newcommand{\cl}{{\cal L}}

\thispagestyle{empty}

\vspace*{1cm}
\begin{center}
{\Large\bf On series of sectorial forms} \\[4mm]

\large C.J.K. Batty$^1$ and A.F.M. ter Elst$^2$

\end{center}

\vspace{4mm}

\begin{center}
{\bf Abstract}
\end{center}

\begin{list}{}{\leftmargin=1.8cm \rightmargin=1.8cm \listparindent=10mm 
   \parsep=0pt}
\item
We prove a convergence theorem for partial sums of sectorial forms
with vertex zero and a common semi-angle.
As an example we prove an absorption theorem for sectorial forms.
\end{list}

\vspace{6mm}
\noindent
June 2013.

\vspace{3mm}

\vspace{3mm}

\vspace{6mm}

\noindent
{\bf Home institutions:}    \\[3mm]
\begin{tabular}{@{}cl@{\hspace{10mm}}cl}
1. & St.\ John's College & 
2. & Department of Mathematics  \\
& University of Oxford & 
  & University of Auckland  \\
& Oxford OX1 3JP  &
  & Private bag 92019 \\
& UK &
   & Auckland 1142 \\ 
&  &
  & New Zealand \\[8mm]
\end{tabular}

\newpage
\setcounter{page}{1}

\section{Introduction} \label{Sabs1}

Kato and Simon proved the following celebrated monotone convergence
theorem for positive symmetric closed forms in a Hilbert space.

\begin{thm} \label{tabs101}
Let $0 \leq \gothh_1 \leq \gothh_2 \leq \ldots$ be positive symmetric closed
sesquilinear forms in a Hilbert space $H$.
Define the form $\gothh_\infty$ by
\[
D(\gothh_\infty)
= \{ u \in \bigcap_{n=1}^\infty D(\gothh_n) : \sup_{n \in \Ni} \gothh_n(u) < \infty \}
\]
and $\gothh_\infty(u,v) = \lim_{n \to \infty} \gothh_n(u,v)$.
Then $\gothh_\infty$ is a closed positive symmetric form.
Suppose $D(\gothh_\infty)$ is dense in $H$.
Let $A_n$ and $A_\infty$ be the self-adjoint operators associated 
with $\gothh_n$ and $\gothh_\infty$ for all $n \in \Ni$.
Then $\lim_{n \to \infty} A_n = A_\infty$ in the strong resolvent sense.
\end{thm}
(See \cite{Kat1} Theorem~VIII.3.13a and \cite{bSim5} Theorem~3.1.)

The condition $D(\gothh_\infty)$ is dense in $H$ in Theorem~\ref{tabs101}
is not essential; without this condition one has to replace
`self-adjoint operator' by `self-adjoint graph' everywhere, see 
\cite{bSim5} Theorem~4.1.
The proofs first use a convergence theorem of Kato, which heavily uses the 
symmetry of the sesquilinear forms.
This convergence theorem states that the sequence of self-adjoint operators
converges in the strong resolvent sense to the resolvent of a self-adjoint 
operator.
In a second step Kato and Simon use different methods to identify this 
limit self-adjoint operator with the operator associated with the form 
$\gothh_\infty$.

The aim of this paper is to extend these results to sectorial forms,
which are possibly not closed, even not closable.
A sesquilinear form $\gota$ in $H$ is called {\bf sectorial}
if there are $\gamma \in \Ri$, called a {\bf vertex},
and $\theta \in [0,\frac{\pi}{2})$, called a {\bf semi-angle}, such that 
\begin{equation}
\gota(u) - \gamma \, \|u\|_H^2
\in \Sigma_\theta
\label{eSabs1;1}
\end{equation}
for all $u \in D(\gota)$, where
\[
\Sigma_\theta
= \{ r \, e^{i \alpha} : r \geq 0, \; |\alpha| \leq \theta \}
.  \]
If $\gota$ is a densely defined sectorial form in $H$, then 
one can associate an $m$-sectorial operator $A$ with $\gota$ in
the following way.
Let $x,f \in H$. 
Then $x \in D(A)$ and $A x = f$ if and only if there exists a 
sequence $(u_n)_{n \in \Ni}$ in $D(\gota)$ such that 
$\lim u_n = x$ in $H$, the set $ \{ \gota(u_n) : n \in \Ni \} $ is 
bounded and for all $v \in D(\gota)$ it follows that 
$\lim \gota(u_n,v) = (f,v)_H$.
(See \cite{AE2} Theorem~1.1.)
We emphasize that the operator $A$ is well-defined (single valued).
Define the norm $\|\cdot\|_{D(\gota)} \colon D(\gota) \to [0,\infty)$ by 
\begin{equation}
\|u\|_{D(\gota)}^2 = \RRe \gota(u) + (1 - \gamma) \|u\|_H^2
\label{eSabs1;10}
\end{equation}
where $\gamma$ is as in (\ref{eSabs1;1}).
The form $\gota$ is called {\bf closed} if the normed space
$(D(\gota),\|\cdot\|_{D(\gota)})$ is complete.

The main theorem of this paper is the following.

\begin{thm} \label{tabs102}
Let $H$ be a Hilbert space.
Fix $\theta \in [0,\frac{\pi}{2})$.
For all $n \in \Ni$ let $\gotb_n$ be a sectorial form in 
$H$ with vertex $0$ and semi-angle~$\theta$.
For all $n \in \Ni$ define 
\[
\gota_n = \sum_{k=1}^n \gotb_k
.  \]
Define 
\[
D(\gota_\infty)
= \{ u \in \bigcap_{n=1}^\infty D(\gotb_n) : \sum \RRe \gotb_n(u) \mbox{ is convergent} \}
.  \]
Then for all $u,v \in D(\gota_\infty)$ the series 
$\sum \gotb_n(u,v)$ is convergent.
Suppose that $D(\gota_\infty)$ is dense in $H$.
Let $A_n$ be the $m$-sectorial operator associated with $\gota_n$
for all $n \in \Ni$.
Then one has the following.
\begin{tabel}
\item \label{tabs102-1}
There exists an $m$-sectorial operator $A_\infty$ in $H$ such that 
$\lim_{n \to \infty} A_n = A_\infty$ in the strong resolvent sense.
\item \label{tabs102-2}
If $\gota_n$ is closed for all $n \in \Ni$ then $\gota_\infty$ is closed.
Moreover, the $m$-sectorial operator $A_\infty$ in Statement~{\rm \ref{tabs102-1}}
is associated with the form $\gota_\infty$.
\end{tabel}
\end{thm}

As in Theorem~\ref{tabs101} the condition $D(\gota_\infty)$ is dense in $H$ is 
not needed if one replaces `operator' by `graph' everywhere in the theorem.
Even if all forms $D(\gota_n)$ are dense in $H$, then $A_\infty$ may be
multi-valued. 
Moreover, if $A_\infty$ is single-valued, then $D(\gota_\infty)$ might not be 
dense in $H$.

Theorem~\ref{tabs102}\ref{tabs102-1} was proved by Ouhabaz \cite{Ouh7} Theorem~5
under the additional condition that either $\IIm \gota_{n+1}(u) \leq \IIm \gota_n(u)$ 
for all $n \in \Ni$ and $u \in D(\gota_{n+1})$
or $\IIm \gota_{n+1}(u) \geq \IIm \gota_n(u)$ 
for all $n \in \Ni$ and $u \in D(\gota_{n+1})$, using Theorem~\ref{tabs101} and 
Vitali's theorem.

The closedness condition in Statement~\ref{tabs102-2} is in general essential. 
In Example~\ref{xabsorp305} we present an example where 
$\gota_n$ is closable for all $n \in \Ni$, the form 
$\gota_\infty$ is closed and densely defined, 
but $A_\infty$ is not associated with the form~$\gota_\infty$.

\medskip

Finally we consider absorption, where we choose $\gotb_n = \gotb_2$ 
for all $n \in \Ni$ with $n \geq 3$.
In that case the density of the form $D(\gota_\infty)$ is usually 
violated, so we phrase the next theorem with a graph instead of operator
for $A_\infty$.

\begin{thm} \label{tabs103}
Let $\gota$ and $\gotb$ be two densely defined sectorial forms in a Hilbert space $\ch$ 
with $D(\gota) = D(\gotb)$ and 
suppose that $\gotb$ has vertex~$0$.
For all $n \in \Ni$ define $\gota_n = \gota + n \, \gotb$
and let $A_n$ be the $m$-sectorial operator associated with $\gota_n$.
Let $A_\infty$ be the $m$-sectorial graph such that 
$\lim A_n = A_\infty$ in the strong resolvent sense.
Further, let $A$ be the $m$-sectorial operator associated with $\gota$.
Then one has the following.
\begin{tabel}
\item \label{tabs103-1.5}
Suppose $\gota$ is closable and there exist $c_1,c_2 > 0$ such that 
\begin{equation}
|\gotb(u)| \leq c_1 \RRe \gota(u) + c_2 \, \|u\|_H^2
\label{etabs103;3}
\end{equation}
for all $u \in \cd$. 
Then there exists an orthogonal projection $P$ in $H$ such that
\[
e^{-t A_\infty}
= \lim_{n \to \infty} \Big( e^{- \frac{t}{n} \, A} \, P \Big)^n
\]
strongly for all $t > 0$.
If in addition $\gota$ is closed, then 
$P$ is the orthogonal projection of $H$ onto 
$ \{ u \in D(\gota) : \gotb(u) = 0 \} \overline{\raisebox{7pt}{$\;\;$}}$, 
where the closure is in $H$.

\item \label{tabs103-2}
Suppose there exists a $($bounded$)$ $B \in \cl(H)$ such that 
\[
\gotb(u,v) = (B u,v)_H
\]
for all $u,v \in D(\gota)$.
Let $P$ be the orthogonal projection from $H$ onto 
$\Big( (B + B^*)(H) \Big)^\perp$.
Then 
\[
e^{-t A_\infty}
= \lim_{n \to \infty} \Big( e^{- \frac{t}{n} \, A} \, P \Big)^n
\]
strongly for all $t > 0$.
\end{tabel}
\end{thm}

\vertspace

Note that in Statement~\ref{tabs103-1.5} it is not assumed that $\gotb$ is closable.
An example is $H = L_2(\Ri)$, $D(\gota) = W^{1,2}(\Ri)$, 
$\gota(u,v) = \int u' \, \overline{v'}$ and $\gotb(u,v) = u(0) \, \overline{v(0)}$.
On the other hand, if $\gota$ is closed and $\gotb$ is closable, then (\ref{etabs103;3})
is valid by \cite{Ouh5} Proposition~1.18.

The results in Statements~\ref{tabs103-1.5} and \ref{tabs103-2}
of Theorem~\ref{tabs103} are optimal.
In Part~\ref{tabs103-1.5} the form $\gotb$ does not have to be bounded 
in $H$, but merely $\gota$-bounded, with the form $\gota$ being closable.
In \ref{tabs103-2} the form $\gota$ does not have to be closable, but the 
form $\gotb$ is bounded in $H$.
In both cases there is a Trotter--Kato-type formula with a projection~$P$.
In Example~\ref{xabsorp102} we present a positive symmetric 
non-closable densely defined form $\gota$ and a positive symmetric 
$\gota$-bounded form $\gotb$ for which there is no
Trotter--Kato-type formula with a projection~$P$.

\medskip

In Section~\ref{Sabs2} we prove Theorem~\ref{tabs102}\ref{tabs102-2}
and in Section~\ref{Sabs3} we prove Theorem~\ref{tabs102}\ref{tabs102-1},
both without the density assumption on the form domain $D(\gota_\infty)$.
In addition we give a short overview about $m$-sectorial graphs in 
Section~\ref{Sabs2}.
In Section~\ref{Sabs4} we prove the absorption in Theorem~\ref{tabs103}.

\section{Sums of closed forms} \label{Sabs2}

In this section we prove Theorem~\ref{tabs102}\ref{tabs102-2}, 
without the density assumption on $D(\gota_\infty)$.
To this end we need the concept of graphs.
For a more thorough introduction we refer to \cite{Bre2}, \cite{Sho}
and \cite{AEKS}.
Fix a Hilbert space $H$.
A {\bf graph} in $H$ is a subspace of $H \times H$.
Let $A$ be a graph.
The {\bf domain} of $A$ is 
$D(A) = \{ x \in H : (\{ x \} \times H) \cap A \neq \emptyset \} $. 
The graph $A$ is called {\bf single-valued} or an {\bf operator}
if $ \{ y \in H : (0,y) \in A \} = \{ 0 \} $.
For operators we will use the usual terminology and notation.
The graph $A$ is called {\bf surjective} if for all $y \in H$ there 
exists an $x \in H$ such that $(x,y) \in A$.
The graph $A$ is called {\bf invertible} if it is 
surjective, closed and the 
reflected graph $ \{ (y,x) : (x,y) \in A \} $ is single-valued.
If the {\em graph} $A$ is invertible then we define the {\em operator}
$A^{-1} \colon H \to H$ by $A^{-1} y = x$ if $(x,y) \in A$.
If $\lambda \in \Ci$ we define the graph
$(A + \lambda \, I)$ by 
\[
(A + \lambda \, I) = \{ (x, y + \lambda \, x) : (x,y) \in A \}
.  \]
The {\bf resolvent} $\rho(A)$ of $A$ is the set of all $\lambda \in \Ci$ such that 
$(A - \lambda \, I)$ is invertible.
The graph $A$ is called {\bf $m$-sectorial} if there are $\gamma \in \Ri$ and 
$\theta \in [0,\frac{\pi}{2})$ such that 
$(x,y)_H \in \Sigma_\theta$ for all $(x,y) \in (A - \gamma \, I)$ and 
$A - (\gamma-1) I$ is invertible.
If $A$ is $m$-sectorial we define the 
{\bf single-valued part $A^\circ$ of $A$} by 
\[
A^\circ = A \cap (\overline{D(A)} \times \overline{D(A)})
.  \]
It is easy to verify that $A^\circ$ is an operator in $\overline{D(A)}$,
it is again $m$-sectorial and $A = A^\circ \oplus ( \{ 0 \} \times D(A)^\perp )$.
Conversely, if $H_1$ is a closed subspace of $H$ and $B$ is an 
$m$-sectorial operator in $H_1$, then $A = B \oplus ( \{ 0 \} \times H_1^\perp )$
is an $m$-sectorial graph in $H$ such that $A^\circ = B$.

Let $\gota$ be a closed sectorial form in $H$.
Define 
\[
A = \{ (u,f) \in D(\gota) \times H : \gota(u,v) = (f,v) \mbox{ for all } v \in D(\gota) \}
.  \]
Then it follows from \cite{Kat1} Theorem~VI.2.1 applied to 
the Hilbert space $\overline{D(\gota)}$ that $A$ is an $m$-sectorial graph.
We call $A$ the {\bf $m$-sectorial graph associated with $\gota$}.

\vertspace

Fix $\theta \in [0,\frac{\pi}{2})$.
For all $n \in \Ni$ let $\gotb_n$ be a sectorial form in 
$H$ with vertex $0$ and semi-angle~$\theta$.
For all $n \in \Ni$ define 
\[
\gota_n = \sum_{k=1}^n \gotb_k
.  \]
Note that $D(\gota_n) = \bigcap_{k=1}^n D(\gotb_k)$ for all $n \in \Ni$.
Next, define 
\begin{equation}
D(\gota_\infty)
= \{ u \in \bigcap_{n=1}^\infty D(\gotb_n) : 
    \sum_{n=1}^\infty \RRe \gotb_n(u) < \infty \} 
. 
\label{eSabs2;1}
\end{equation}
If $u,v \in D(\gota_\infty)$ then 
\[
|\gotb_n(u,v)| 
\leq (1 + \tan \theta) (\RRe \gotb_n(u))^{1/2} \, (\RRe \gotb_n(v))^{1/2}
\leq (1 + \tan \theta) \Big( \RRe \gotb_n(u) + \RRe \gotb_n(v) \Big)
\]
for all $n \in \Ni$.
So $\sum \gotb_n(u,v)$ is convergent.
Define $\gota_\infty \colon D(\gota_\infty) \times D(\gota_\infty) \to \Ci$ by 
\[
\gota_\infty(u,v) = \sum_{n=1}^\infty \gotb_n(u,v)
.  \]
Then 
$|\IIm \gota_\infty(u)| \leq (\tan \theta) \RRe \gota_\infty(u)$
for all $u \in D(\gota_\infty)$.
So $\gota_\infty$ is sectorial with vertex $0$ and semi-angle~$\theta$.

\begin{lemma} \label{labsorp201}
Suppose the form $\gota_n$ is closed for all $n \in \Ni$.
Then $\gota_\infty$ is closed.
\end{lemma}
\proof\
This follows from \cite{bSim5} Theorem~4.1.
For self-consistency we give a proof.
Let $(u_n)_{n \in \Ni}$ be a Cauchy sequence in $D(\gota_\infty)$.
Then $(u_n)_{n \in \Ni}$ is a Cauchy sequence in $H$, so 
$u = \lim_{n \to \infty} u_n$ exists in $H$.
Let $m \in \Ni$. 
Then $(u_n)_{n \in \Ni}$ is a Cauchy sequence in $D(\gota_m)$ and 
$\gota_m$ is closed. 
Hence $u \in D(\gota_m)$.
Moreover, $\RRe \gota_m(u) \leq \liminf_{n \to \infty} \RRe \gota_m(u_n)
\leq \liminf_{n \to \infty} \RRe \gota_\infty(u_n)$.
So $u \in \bigcap_{m=1}^\infty D(\gota_m)$ and 
$\sup_{m \in \Ni} \RRe \gota_m(u) \leq \liminf_{n \to \infty} \RRe \gota_\infty(u_n)$.
Therefore $u \in D(\gota_\infty)$.
Finally, let $n \in \Ni$.
Then 
\[
\RRe \gota_m(u_n - u)
\leq \liminf_{k \to \infty} \RRe \gota_m(u_n - u_k)
\leq \liminf_{k \to \infty} \RRe \gota_\infty(u_n - u_k)
\]
for all $m \in \Ni$.
So 
\[
\RRe \gota_\infty(u_n - u) \leq \sup_{m \in \Ni} \RRe \gota_m(u_n - u)
\leq \liminf_{k \to \infty} \RRe \gota_\infty(u_n - u_k)
.  \]
Therefore $\lim_{n \to \infty} \RRe \gota_\infty(u_n - u) = 0$
and $D(\gota_\infty)$ is closed.\hfill$\Box$

\vertspace

Theorem~\ref{tabs102}\ref{tabs102-2} is a special case of 
Lemma~\ref{labsorp201} and the next theorem.

\begin{thm} \label{tabsorp202}
Let $H$ be a Hilbert space and $\theta \in [0,\frac{\pi}{2})$.
For all $n \in \Ni$ let $\gotb_n$ be a sectorial form in 
$H$ with vertex $0$ and semi-angle~$\theta$.
For all $n \in \Ni$ define 
$\gota_n = \sum_{k=1}^n \gotb_k$.
Define $D(\gota_\infty)$ as in {\rm (\ref{eSabs2;1})} and 
define $\gota_\infty = \sum_{n=1}^\infty \gotb_n$.
Suppose that $\gota_n$ is closed for all $n \in \Ni$.
For all $n \in \Ni$ let $A_n$ be the graph associated with $\gota_n$
and let $A_\infty$ be the graph associated with~$\gota_\infty$.
Then $\lim A_n = A_\infty$ in the strong resolvent sense.
\end{thm}
\proof\
Without loss of generality we may assume that 
$D(\gotb_{n+1}) \subset D(\gotb_n)$ for all $n \in \Ni$.

Let $f \in H$.
For all $n \in \Ni$ set $u_n = (A_n + I)^{-1} f$.
Then $u_n \in D(\gota_n)$ and 
\begin{equation}
\gota_n(u_n,v) + (u_n,v)_H = (f,v)_H
\label{etabsorp202;2}
\end{equation}
for all $v \in D(\gota_n)$.
Choose $v = u_n$.
Then 
\[
\RRe \gota_n(u_n) + \|u_n\|_H^2 = \RRe (f,u_n)_H \leq \|f\|_H \, \|u_n\|_H
.  \]
So $\|u_n\|_H \leq \|f\|_H$ and $\RRe \gota_n(u_n) \leq \|f\|_H^2$.
Let $m \in \Ni$.
If $n \in \Ni$ and $m \leq n$, then 
$\RRe \gota_m(u_n) \leq \RRe \gota_n(u_n) \leq \|f\|_H^2$.
So the (tail of the) sequence $(u_n)_{n \in \Ni}$ is bounded in $D(\gota_m)$.
Using a diagonal argument and passing to a subsequence if necessary, 
we may assume that for all $m \in \Ni$ the sequence $(u_n)_{n \in \Ni}$
is weakly convergent in $D(\gota_m)$.
Since $D(\gota_m)$ is continuously embedded in the Hausdorff space $H$,
there exists a $u \in \bigcap_{m=1}^\infty D(\gota_m)$ such that 
$(u_n)_{n \in \Ni}$ converges weakly to $u$ in $D(\gota_m)$ for all $m \in \Ni$.
Then 
\begin{equation}
\RRe \gota_m(u) 
\leq \liminf_{n \to \infty} \RRe \gota_m(u_n) 
\leq \liminf_{n \to \infty} \RRe \gota_n(u_n) 
\leq \|f\|_H^2
\label{etabsorp202;1}
\end{equation}
for all $m \in \Ni$.
Therefore 
$\sum_{n=1}^\infty \RRe \gotb_n(u)
= \sup_{m \in \Ni} \RRe \gota_m(u) \leq \|f\|_H^2$
and $u \in D(\gota_\infty)$.

Let $v \in D(\gota_\infty)$.
Then $\gota_n(u_n,v) + (u_n,v)_H = (f,v)_H$ for all $n \in \Ni$ by (\ref{etabsorp202;2}).
We shall show that $\lim \gota_n(u_n,v) = \gota_\infty(u,v)$.
Clearly $\lim \gota_n(u,v) = \gota_\infty(u,v)$, so it 
suffices to show that $\lim \gota_n(u_n-u,v) = 0$.
Let $\varepsilon > 0$.
Because $v \in D(\gota_\infty)$ there exists an $N_0 \in \Ni$ such that 
$\sum_{k=N_0}^\infty \RRe \gotb_k(v) \leq \varepsilon^2$.
Since $\lim_{n \to \infty} \gota_{N_0}(u_n - u,v) = 0$, there exists an 
$N \geq N_0$ such that $|\gota_{N_0}(u_n - u,v)| \leq \varepsilon$
for all $n \geq N$.
Then for all $n \geq N$ one estimates
\begin{eqnarray*}
|\gota_n(u_n-u,v)|
& \leq & |\gota_{N_0}(u_n - u,v)| + \sum_{k=N_0+1}^n |\gotb_k(u_n-u,v)|  \\
& \leq & \varepsilon 
    + \sum_{k=N_0+1}^n (1 + \tan \theta) (\RRe \gotb_k(u_n-u))^{1/2} \, (\RRe \gotb_k(v))^{1/2}  \\
& \leq & \varepsilon 
    + \sum_{k=N_0+1}^n (1 + \tan \theta) 
         \Big( \varepsilon \RRe \gotb_k(u_n-u) + \varepsilon^{-1} \, \RRe \gotb_k(v) \Big) \\
& \leq & \varepsilon 
    +  (1 + \tan \theta) 
         \Big( \varepsilon \RRe \gota_n(u_n-u) 
                 + \varepsilon^{-1} \, \sum_{k=N_0+1}^\infty \RRe \gotb_k(v) \Big)  \\
& \leq & \varepsilon 
    +  (1 + \tan \theta) 
         \Big( \varepsilon (2 \RRe \gota_n(u_n) + 2 \RRe \gota_n(u) ) 
                 + \varepsilon \Big)  \\
& \leq & (1 + (1 + \tan \theta) (4 \|f\|_H^2 + 1)) \, \varepsilon
.
\end{eqnarray*}
So $\gota_\infty(u,v) + (u,v)_H = (f,v)_H$.
Hence $(u,f) \in A_\infty + I$ and $(A_\infty + I)^{-1} f = u$.
Since $\lim u_n = u$ weakly in $H$, one has $\|u\|_H \leq \liminf \|u_n\|_H^2$.
In order to prove strong convergence it remains to show that 
$\limsup \|u_n\|_H^2 \leq \|u\|_H$.
It follows from (\ref{etabsorp202;1}) that 
$\RRe \gota_\infty(u) = \sup_{m \in \Ni} \RRe \gota_m(u)
\leq \liminf_{n \to \infty} \RRe \gota_n(u_n)$.
So (\ref{etabsorp202;2}) gives
\begin{eqnarray*}
\limsup_{n \to \infty} \|u_n\|_H^2
& = & \limsup_{n \to \infty} \Big( \RRe (f,u_n)_H - \RRe \gota_n(u_n) \Big)
\leq \RRe (f,u)_H - \RRe \gota_\infty(u)
= \|u\|_H^2
.
\end{eqnarray*}
Therefore $\lim u_n = u$ strongly.\hfill$\Box$

\section{Sums of non-closed forms} \label{Sabs3}

In this section we prove Theorem~\ref{tabs102}\ref{tabs102-1}, 
without the density assumption on $D(\gota_\infty)$.
For the statement and its proof we need to extend two definitions
and theorems of \cite{AE2} to the setting of $m$-sectorial 
graphs instead of $m$-sectorial operators.

Let $\gota$ be a sectorial sesquilinear form in a Hilbert space $H$.
The {\bf graph associated with $\gota$} is the set of all $(x,f) \in H \times H$
for which there exists a 
sequence $(u_n)_{n \in \Ni}$ in $D(\gota)$ such that 
$\lim u_n = x$ in $H$, the set $ \{ \gota(u_n) : n \in \Ni \} $ is 
bounded and for all $v \in D(\gota)$ it follows that 
$\lim \gota(u_n,v) = (f,v)_H$.
Then it follows from \cite{AE2} Theorem~1.1 applied to $\overline{D(\gota)}$
that $A$ is an $m$-sectorial graph.

Next we wish to consider sesquilinear forms for which the 
form domain is no longer a subspace of $H$.
Let $V$ and $H$ be two Hilbert spaces and $\gota \colon V \times V \to \Ci$ 
be a sesquilinear form. 
The form $\gota$ is called {\bf continuous} if there exists an $M > 0$ 
such that $|\gota(u,v)| \leq M \, \|u\|_V \, \|v\|_V$
for all $u,v \in V$.
Let $j \colon V \to H$ be a continuous linear map.
The form $\gota$ is called {\bf $j$-elliptic} if 
there exist $\omega \in \Ri$ and $\mu > 0$ such that 
\[
\RRe \gota(u) + \omega \, \|j(u)\|_H^2
\geq \mu \, \|u\|_V^2
\]
for all $u \in V$.
The {\bf graph associated with $(\gota,j)$} is the space
\begin{eqnarray*}
A = \{ (x,f) \in H \times H & : & 
          \mbox{there exists a $u \in V$ such that}    \\
& & j(u) = x \mbox{ and } \gota(u,v) = (f,v) 
     \mbox{ for all } v \in V \}
.  
\end{eqnarray*}
If $\gota$ is continuous and $j$-elliptic, then it follows from 
\cite{AE2} Theorem~2.1 applied to the Hilbert space $\overline{j(V)}$
that the graph associated with $(\gota,j)$ is an $m$-sectorial graph.

The following construction is useful to deal with possibly 
non-closed forms.

\begin{lemma} \label{labs310}
Let $\gota$ be a sectorial form in a Hilbert space $H$
with vertex $\gamma$.
Define the norm $\|\cdot\|_{D(\gota)}$ on $D(\gota)$ as in 
{\rm (\ref{eSabs1;10})}.
Let $V$ be a Hilbert space, $j \in \cl(V,H)$, 
$q \colon D(\gota) \to V$ linear and $\tilde \gota \colon V \times V \to \Ci$ 
a continuous $j$-elliptic sesquilinear form.
Suppose $q(D(\gota))$ is dense in $V$ and there exists a $c > 0$ 
such that 
\[
c^{-1} \, \|q(u)\|_V \leq \|u\|_{D(\gota)} \leq c \, \|q(u)\|_V
,  \]
$j(q(u)) = u$ and $\tilde \gota(q(u), q(v)) = \gota(u,v)$
for all $u,v \in D(\gota)$.

Then the graph associated with $\gota$ is equal to the graph
associated with $(\tilde \gota,j)$.
\end{lemma}
\proof\
This follows from \cite{AE2} Proposition~3.3 applied to 
$\overline{D(\gota)} = \overline{j(V)}$.\hfill$\Box$

\vertspace

Now we are able to state and prove the alluded extension of 
Theorem~\ref{tabs102}\ref{tabs102-1} to forms which are possibly 
not densely defined.

\begin{thm} \label{tabsorp301}
Let $H$ be a Hilbert space.
Fix $\theta \in [0,\frac{\pi}{2})$.
For all $n \in \Ni$ let $\gotb_n$ be a sectorial form in 
$H$ with vertex $0$ and semi-angle~$\theta$.
For all $n \in \Ni$ define 
\[
\gota_n = \sum_{k=1}^n \gotb_k
.  \]
Let $A_n$ be the $m$-sectorial graph associated with $\gota_n$
for all $n \in \Ni$.
Then there exists an $m$-sectorial graph $A_\infty$ in $H$ such that 
$\lim_{n \to \infty} A_n = A_\infty$ in the strong resolvent sense.
\end{thm}
\proof\
The proof of the theorem is a modification of the 
proof of Theorem~\ref{tabsorp202}, with several complications.
For each $n$ one can consider the completion of 
$(D(\gota_n), \|\cdot\|_{D(\gota_n)})$.
Unfortunately the natural maps between these completions are 
in general not 
injective, whilst in the proof of Theorem~\ref{tabsorp202} the
spaces $D(\gota_n)$ are decreasing.
The operator $A_\infty$ is associated with a form on a Hilbert 
space and for this Hilbert space we cannot take the direct sum 
of the completions.
Therefore we put a weighted norm on $D(\gota_n)$ such that 
in the end the graph $A_\infty$ can be associated with a 
form on a closed subspace of the direct sum, together with a 
continuous map~$j$.
We lift the forms $\gota_n$ and $\gotb_k$ to the level of the completions
and denote these forms with tildes.
Then we lift them to the direct sum and denote them with hats.

Without loss of generality we may assume that 
$D(\gotb_{n+1}) \subset D(\gotb_n)$ for all $n \in \Ni$.

Define $\gota_0 = \gotb_0 \colon H \times H \to \Ci$ by 
$\gotb_0(u,v) = (u,v)_H$.
Let $n \in \Ni_0$.
Define the norm $\|\cdot\|_n \colon D(\gota_n) \to [0,\infty)$ by 
\[
\|u\|_n^2 = \sum_{k=0}^n 2^{-(n-k)} \, \RRe \gotb_k(u)
.  \]
Let $V_n$ be the completion of $(D(\gota_n), \|\cdot\|_n)$ and 
let $q_n \colon D(\gota_n) \to V_n$ be the natural map.
Note that $V_0 = H$, with the same norm.

Let $m,n \in \Ni_0$ with $m \leq n$.
Then the inclusion $i_{nm} \colon (D(\gota_n), \|\cdot\|_n) \to (D(\gota_m), \|\cdot\|_m)$
is continuous. 
Hence there exists a unique continuous $\Phi_{nm} \colon V_n \to V_m$
such that $\Phi_{nm} \circ q_n = q_m \circ i_{nm}$.
Note that $\Phi_{n0} \in \cl(V_n,H)$ and $\Phi_{n0} (q_n(u)) = u$
for all $u \in D(\gota_n)$.

Let $k,n \in \Ni_0$ with $k \leq n$.
There exists a unique continuous sesquilinear form 
$\tilde \gotb_{nk} \colon V_n \times V_n \to \Ci$
such that 
\[
\tilde \gotb_{nk}(q_n(u), q_n(v)) = \gotb_k(u,v)
.   \]
Then $\tilde \gotb_{nk}$ is sectorial with vertex $0$ and semi-angle~$\theta$.
Note that if $m \in \Ni_0$ is such that $k \leq m \leq n$, then 
$\tilde \gotb_{nk}(u) = \tilde \gotb_{mk}(\Phi_{nm}(u))$
for all $u \in V_n$.
Define $\tilde \gota_n \colon V_n \times V_n \to \Ci$ by 
\[
\tilde \gota_n = \sum_{k=1}^n \tilde \gotb_{nk}
.  \]
It is easy to verify that 
\begin{equation}
\|u\|_{V_n}^2
= \sum_{k=0}^n 2^{-(n-k)} \, \RRe \tilde \gotb_{nk}(u)
\label{etabsorp301;40}
\end{equation}
for all $u \in V_n$.

Next, define $V = \oplus_{n=0}^\infty V_n$.
For all $n \in \Ni_0$ let $\pi_n \colon V \to V_n$ be the natural projection.
Define $j \in \cl(V,H)$ by $j = \pi_0$.
Let $n \in \Ni_0$.
Define the closed space
\begin{eqnarray*}
W_n = \{ (u^{(\ell)})_{\ell \in \Ni_0} \in V 
& : &  
           u^{(\ell)} = 0 \mbox{ for all } \ell \in \{ n+1,n+2,\ldots \}  \mbox{ and }   \\*
& &            u^{(\ell)} = \Phi_{n\ell}(u^{(n)}) \mbox{ for all } \ell \in \{ 0,\ldots,n-1 \} \}
.  
\end{eqnarray*}
Define $I_n \colon V_n \to W_n$ by 
\[
I_n(u) = (\Phi_{n0}(u), \ldots,\Phi_{nn}(u),0,0,\ldots)
.  \]
Then $I_n$ is a continuous bijection.
Moreover, $(j \circ (I_n \circ q_n))(u) = u$ for all $u \in D(\gota_n)$
and  $(I_n \circ q_n)(D(\gota_n))$ is  dense in $W_n$.
For all $k \in \Ni_0$ define $\hat \gotb_k \colon V \times V \to \Ci$ by 
\[
\hat \gotb_k(u,v)
= \tilde \gotb_{kk}(\pi_k(u),\pi_k(v))
.  \]
Then $\hat \gotb_k$ is sectorial with vertex $0$ and semi-angle~$\theta$.
For all $n \in \Ni$ define $\hat \gota_n \colon V \times V \to \Ci$ by
\[
\hat \gota_n = \sum_{k=1}^n \hat \gotb_k
\]
and define $\hat \gota_0 \colon V \times V \to \Ci$ by $\hat \gota_0 = \hat \gotb_0$.
If $n \in \Ni$ then 
\begin{equation}
\hat \gota_n(u) 
= \sum_{k=1}^n \hat \gotb_k(u)
= \sum_{k=1}^n \tilde \gotb_{kk}(\pi_k(u))
= \sum_{k=1}^n \tilde \gotb_{nk}(\pi_n(u))
= \tilde \gota_n(\pi_n(u))
\label{etabsorp301;4}
\end{equation}
for all $u \in W_n$.
Moreover, $\hat \gota_n((I_n \circ q_n)(u)) = \gota_n(u)$
for all $u \in D(\gota_n)$.
It is easy to see that 
\[
\|u\|_{W_n}^2
= \sum_{m=0}^n \sum_{k=0}^m 2^{-(m-k)} \, \RRe \hat \gotb_m(u)
\]
for all $u \in W_n$.
Therefore the form $\hat \gota_n|_{W_n \times W_n}$ is continuous and 
$j|_{W_n}$-elliptic.
Moreover, $A_n$ is the graph associated with 
$(\hat \gota_n|_{W_n \times W_n}, j|_{W_n})$ by Lemma~\ref{labs310}.

For later purposes, for all $m \in \Ni_0$ define the truncation
$T_m \colon V \to V$ by
\[
T_m( u^{(0)}, u^{(1)},\ldots ) 
= (u^{(0)}, \ldots, u^{(m)}, 0, 0,\ldots )
.  \]
Then $T_m$ is continuous.
Moreover, if $n \in \{ m,m+1,\ldots \} $ then 
$T_m(W_n) \subset W_m$ and $\hat \gota_m(T_m u) = \hat \gota_m(u)$
for all $u \in W_n$.

Define 
\[
W_\infty
= \{ (u^{(\ell)})_{\ell \in \Ni_0} \in V : 
        \Phi_{nm}(u^{(n)}) = u^{(m)} \mbox{ for all } n,m \in \Ni_0
        \mbox{ with } m \leq n \}
.  \]
We need a lemma.

\begin{lemma} \label{labsorp302}
For all $n \in \Ni_0$ let $u_n \in V_n$.
Then the following are equivalent.
\begin{tabeleq}
\item \label{labsorp302-1}
$(u_0,u_1,\ldots) \in W_\infty$.
\item \label{labsorp302-2}
$\Phi_{nm}(u_n) = u_m$ for all $n,m \in \Ni_0$ with $m \leq n$ and 
$\{ \RRe \tilde \gota_n(u_n) : n \in \Ni \} $ is bounded.
\end{tabeleq}
Moreover, if {\rm \ref{labsorp302-1}} is valid, then 
\begin{equation}
\|(u_0,u_1,\ldots)\|_V^2 \leq 2 (\sup_{n \in \Ni} \RRe \tilde \gota_n(u_n) + \|u_0\|_H^2)
.
\label{elabsorp302;1}
\end{equation}
\end{lemma}
\proof\
`\ref{labsorp302-1}$\Rightarrow$\ref{labsorp302-2}'.
Let $n \in \Ni$.
Then $\RRe \tilde \gotb_{nn}(u_n) \leq \|u_n\|_{V_n}^2$.
So if $k \in \{ 1,\ldots,n \} $ then
$\tilde \gotb_{nk}(u_n) = \tilde \gotb_{kk}(\Phi_{nk}(u_n)) = \tilde \gotb_{kk}(u_k)$.
Therefore 
$\RRe \tilde \gota_n(u_n) 
= \sum_{k=1}^n \RRe \tilde \gotb_{nk}(u_n)
\leq \sum_{k=1}^n \|u_k\|_{V_k}^2
\leq \|(u_0,u_1,\ldots)\|_V^2$.

`\ref{labsorp302-2}$\Rightarrow$\ref{labsorp302-1}'.
Let $M = \sup \{ \RRe \tilde \gota_n(u_n) : n \in \Ni \} $.
Note that $\RRe \tilde \gotb_{nk}(u_n) = \RRe \tilde \gotb_{kk}(u_k)$
for all $k,n \in \Ni_0$ with $k \leq n$.
So if $N \in \Ni$ then 
\begin{eqnarray*}
\sum_{n=0}^N \|u_n\|_{V_n}^2
& = & \sum_{n=0}^N \sum_{k=0}^n 2^{-(n-k)} \, \RRe \tilde \gotb_{nk}(u_n)  \\
& = & \sum_{n=0}^N \sum_{k=0}^n 2^{-(n-k)} \, \RRe \tilde \gotb_{kk}(u_k)  \\
& \leq & \sum_{k=0}^N \sum_{j=0}^N 2^{-j} \, \RRe \tilde \gotb_{kk}(u_k)  \\
& \leq & 2 \sum_{k=0}^N \RRe \tilde \gotb_{Nk}(u_N)  \\
& = & 2 (\RRe \tilde \gota_N(u_N) + \|u_0\|_H^2)
\leq 2 (M + \|u_0\|_H^2)
.
\end{eqnarray*}
So $(u_0,u_1,\ldots) \in V$.\hfill$\Box$

\vertspace

It follows from Lemma~\ref{labsorp302} and (\ref{etabsorp301;4}) 
that $\sum \RRe \hat \gotb_n(u)$
is convergent for all $u \in W_\infty$.
If $u,v \in W_\infty$ then 
\[
|\hat \gotb_n(u,v)| 
\leq (1 + \tan \theta) (\RRe \hat \gotb_n(u))^{1/2} \, (\RRe \hat \gotb_n(v))^{1/2}
\leq (1 + \tan \theta) \Big( \RRe \hat \gotb_n(u) + \RRe \hat \gotb_n(v) \Big)
\]
for all $n \in \Ni$.
So $\sum \hat \gotb_n(u,v)$ is convergent.
Define $\hat \gota_\infty \colon W_\infty \times W_\infty \to \Ci$ by 
\[
\hat \gota_\infty(u,v) = \sum_{n=1}^\infty \hat \gotb_n(u,v)
.  \]
Then 
$|\IIm \hat \gota_\infty(u)| \leq (\tan \theta) \RRe \hat \gota_\infty(u)$
for all $u \in W_\infty$.
So $\hat \gota_\infty$ is sectorial with vertex $0$ and semi-angle~$\theta$.
If $u \in W_\infty$ then 
\[
\RRe \hat \gota_\infty(u)
= \sum_{n=1}^\infty \RRe \hat \gotb_n(u)
= \sum_{n=1}^\infty \RRe \tilde \gotb_{nn}(\pi_n(u))
\leq \sum_{n=1}^\infty \|\pi_n(u)\|_{V_n}^2
\leq \|u\|_{W_\infty}^2
.  \]
So $\hat \gota_\infty$ is continuous.
It follows from (\ref{elabsorp302;1}) that $\hat \gota_\infty$ is 
$j|_{W_\infty}$-elliptic.
Let $A_\infty$ be the $m$-sectorial graph associated with 
$(\hat \gota_\infty, j|_{W_\infty})$.
We shall show that the sequence $(A_n)_{n \in \Ni}$ converges in 
the strong resolvent sense to $A_\infty$.

Let $f \in H$.
For all $n \in \Ni$ set $x_n = (A_n + I)^{-1} f$.
Since $A_n$ is the graph associated with $(\hat \gota_n|_{W_n \times W_n}, j|_{W_n})$,
there exists a $u_n \in W_n$ such that $j(u_n) = x_n$ and 
\[
\hat \gota_n(u_n,v) + (j(u_n),j(v))_H = (f,j(v))_H
\]
for all $v \in W_n$.
Choose $v = u_n$.
Then 
\[
\RRe \hat \gota_n(u_n) + \|j(u_n)\|_H^2 = \RRe (f,j(u_n))_H \leq \|f\|_H \, \|j(u_n)\|_H
.  \]
So $\|j(u_n)\|_H \leq \|f\|_H$ and $\RRe \hat \gota_n(u_n) \leq \|f\|_H^2$.

Let $m \in \Ni_0$.
If $n \in \Ni$ and $1 \leq m \leq n$ then 
\[
\RRe \hat \gota_m(T_m u_n)
= \RRe \hat \gota_m(u_n)
\leq \RRe \hat \gota_n(u_n)
\leq \|f\|_H^2
.  \]
Alternatively, if $n \in \Ni$ and $0 = m \leq n$ then 
$\RRe \hat \gota_m(T_m u_n)
= \|j(u_n)\|_H^2
\leq \|f\|_H^2$.
So the (tail of the) sequence $(T_m u_n)_{n \in \Ni}$ is bounded in $W_m$.
Using a diagonal argument and passing to a subsequence if necessary, 
we may assume that for all $m \in \Ni_0$ the sequence $(T_m u_n)_{n \in \Ni}$
is weakly convergent in $W_m$.
Let $m \in \Ni_0$.
Let $w_m \in W_m$ be such that $\lim_{n \to \infty} T_m u_n = w_m$
weakly in $W_m$.
Because $T_m$ is continuous and $T_n \, T_m = T_m$ for all $n \in \Ni_0$ 
with $n \geq m$ one obtains that $T_m w_n = w_m$.
So $\pi_k w_n = \pi_k w_m$ for all $k,m,n \in \Ni_0$ with $k \leq n \wedge m$.
Define $\mu^{(n)} = \pi_n w_n \in V_n$ for all $n \in \Ni_0$.
Then $\Phi_{nm}(\mu^{(n)}) = \mu^{(m)}$ for all $n,m \in \Ni_0$ with $m \leq n$.
Moreover, $I_n(\mu^{(n)}) = w_n$.
So $\RRe \tilde \gota_n(\mu^{(n)}) = \RRe \hat \gota_n(w_n)$
for all $n \in \Ni_0$.
Next, if $m \in \Ni$ then 
$\hat \gota_m|_{W_m \times W_m}$ is continuous and one deduces that 
\[
\RRe \hat \gota_m(w_m)
\leq \liminf_{n \to \infty} \RRe \hat \gota_m(T_m u_n)
\leq \RRe \hat \gota_n(u_n)
\leq \|f\|_H^2
.
\]
In particular, $\RRe \tilde \gota_m(\mu^{(m)}) \leq \|f\|_H^2$.
Therefore 
\[
u = (\mu^{(0)},\mu^{(1)},\ldots )
\in W_\infty
\]
by Lemma~\ref{labsorp302}.
Let $m \in \Ni$.
Then $T_m u = w_m$.
Let $v \in W_m$.
Then $\hat \gota_m(u_n - u,v) = \hat \gota_m(T_m(u_n - u), v) = \hat \gota_m(T_m u_n - w_m, v)$
for all $n \geq m$.
Therefore 
\[
\lim_{n \to \infty} \hat \gota_m(u_n - u,v) = 0
\]
by the weak convergence on $W_m$.
Moreover, the weak convergence on $W_0$ implies that 
$\lim_{n \to \infty} T_0 u_n = w_0$ weakly.
Hence $\lim_{n \to \infty} j(u_n) = \lim_{n \to \infty} j(T_0 u_n) = j(w_0) = j(u)$
weakly in $H$.
Now the rest of the proof is similar to the proof of 
Theorem~\ref{tabsorp202}, with obvious changes.
We leave the details to the reader.\hfill$\Box$

\begin{prop} \label{pabsorp303}
Assume the conditions and notation as in Theorem~{\rm \ref{tabsorp301}}.
Define 
\[
D(\gota_\infty)
= \{ u \in \bigcap_{n=1}^\infty D(\gotb_n) : \sum \RRe \gotb_n(u) \mbox{ is convergent} \}
.  \]
Then $\overline{D(\gota_\infty)} \subset \overline{D(A_\infty)}$.
In particular, if $D(\gota_\infty)$ is dense in $H$ then $A_\infty$ is 
an $m$-sectorial operator. 
\end{prop}
\proof\
We use the notation as in the proof of Theorem~{\rm \ref{tabsorp301}}.
If $u \in D(\gota_\infty)$, then $(q_0(u), q_1(u),\ldots) \in W_\infty$
by Lemma~\ref{labsorp302}.
Moreover, $j((q_0(u), q_1(u),\ldots)) = u$.
Hence $D(\gota_\infty) \subset j(W_\infty)$.
Now the statement follows from \cite{AE2} Theorem~2.5(ii)
and \cite{Kat1} Theorem~VI.2.1.\hfill$\Box$

\vertspace

One can define the form $\gota_\infty$ by 
\[
\gota_\infty(u,v) = \sum_{k=1}^\infty \gota_k(u,v)
\]
for all $u,v \in D(\gota_\infty)$, where we continue to use the 
notation as above.
Even if $D(\gota_\infty)$ is dense in $H$, in general the 
graph/operator $A_\infty$ is {\em not} associated with the form 
$\gota_\infty$.
For an example we first need a lemma.

\begin{lemma} \label{labsorp304}
Let $H$ be an infinite dimensional Hilbert space.
Then there exist dense subspaces $\cd_1,\cd_2,\ldots$ in $H$ 
such that for all $m \in \Ni$ and $u_1 \in \cd_1$,\ldots, $u_m \in \cd_m$
with $u_1 + \ldots + u_m = 0$ it follows that $u_1 = \ldots = u_m = 0$.
\end{lemma}
\proof\
Without loss of generality we may assume that $H$ is separable 
and $H = L_2(0,2\pi)$.
The construction is a modification of the proof of Corollary~1 on page~274
in \cite{FW}.
For all $k \in \Zi$ define $e_k \colon [0,2\pi] \to \Ci$ by 
$e_k(x) = e^{i k x}$ and set 
\[
\cd_0 
= \spann \{ e_k : k \in \Zi \} 
.  \]
Then $\cd_0$ is dense in $H$ and every element of $\cd_0$ can be extended
to an entire function.
For all $m \in \Ni$ define $\varphi_m \colon [0,2\pi] \to \{ -1,1 \} $ by
\[
\varphi_m(x) 
= \left\{ \begin{array}{ll}
   1 & \mbox{if } x \in [0,\tfrac{1}{m}],  \\[5pt]
   -1 & \mbox{if } x \in (\tfrac{1}{m},2\pi] .
          \end{array} \right.
\]
Then $u \mapsto \varphi_m \, u$ is a unitary map in $H$.
Define 
\[
\cd_m = \{ \varphi_m \, u : u \in \cd_0 \}
.  \]
Then $\cd_m$ is dense in $H$.
Finally, let $m \in \Ni$ and $u_1,\ldots,u_m \in \cd_0$.
Suppose that $\sum_{j=1}^m \varphi_j \, u_j = 0$.
Without loss of generality we may assume that $m \geq 2$.
Then 
\[
(\varphi_m \, u_m)|_{(0,\frac{1}{m-1})} 
= - \sum_{j=1}^{m-1} (\varphi_j \, u_j)|_{(0,\frac{1}{m-1})} 
= - \sum_{j=1}^{m-1} u_j|_{(0,\frac{1}{m-1})} 
\]
is the restriction of an entire function to the interval $(0,\frac{1}{m-1})$.
Hence $u_m = 0$.
By induction $u_1 = \ldots = u_{m-1} = 0$.\hfill$\Box$

\begin{exam} \label{xabsorp305}
Let $\Omega$ be the open unit ball in $\Ri^2$.
Let $K = W^{1,2}_0(\Omega)^\perp$, where the orthogonal complement
is in $W^{1,2}(\Omega)$.
Then with the induced inner product of $W^{1,2}(\Omega)$ the space
$K$ is an infinite dimensional Hilbert space.
By Lemma~\ref{labsorp304} there exist 
dense subspaces $\cd_1,\cd_2,\ldots$ in $K$ 
such that for all $m \in \Ni$ and $u_1 \in \cd_1$, \ldots, $u_m \in \cd_m$
with $u_1 + \ldots + u_m = 0$ it follows that $u_1 = \ldots = u_m = 0$.
For all $N \in \Ni$ let 
\[
D(\gotb_N) 
= W^{1,2}_0(\Omega) \oplus \spann \Big( \bigcup_{j=N}^\infty \cd_j \Big)
,  \]
where we emphasize that the span is the finite linear span.
Define 
\[
\gotb_1(u,v)
= \int_\Omega \nabla u \cdot \overline{\nabla v}
\]
for all $u,v \in D(\gotb_1)$.
For all $N \geq 2$ define $\gotb_N(u,v) = 0$ for all $u,v \in D(\gotb_N)$.
Finally, choose $H = L_2(\Omega)$.
Then $\gotb_N$ is a densely defined positive symmetric form in $H$
for all $N \in \Ni$.
We use the notation as in Theorem~{\rm \ref{tabsorp301}}
and Proposition~\ref{pabsorp303}.
Then $D(\gota_\infty) = W^{1,2}_0(\Omega)$
and $\gota_\infty = \gotb_1|_{W^{1,2}_0(\Omega) \times W^{1,2}_0(\Omega)}$.
So $\gota_\infty$ is closed and densely defined.
The operator associated with $\gota_\infty$ is minus the 
Dirichlet Laplacian.

Let $N \in \Ni$.
Then $\gota_N = \gotb_1|_{D(\gotb_N) \times D(\gotb_N)}$.
Since $D(\gotb_N)$ is dense in $W^{1,2}(\Omega)$, the
operator $A_N$ associated with $\gota_N$ is minus the 
Neumann Laplacian.

It is trivial to see that the operators $A_N$ converge to 
minus the Neumann Laplacian in the strong resolvent sense. 
This operator is not equal to the operator associated with $\gota_\infty$.
\end{exam}

This example is a big contrast to Theorem~3.7 in \cite{AE2},
where `convergence from above' is proved for sectorial forms which do not have 
to be closable and one always has convergence in the strong resolvent sense.

\section{Absorption} \label{Sabs4}

In this section we prove Theorem~\ref{tabs103}  and extend it 
to graphs.
Without loss of generality we may assume that $\gota$ has vertex $0$ too.

Let $A$ be an $m$-sectorial graph in a Hilbert space $H$.
Then for all $t > 0$ define the operator $e^{-t A}$ by 
\[
e^{-tA} = \lim_{n \to \infty} \Big( (I + \tfrac{t}{n} \, A)^{-1} \Big)^n
.  \]
Let $A^\circ$ be the single-valued part of $A$.  
Using the decomposition $H = \overline{D(A)} \oplus D(A)^\perp$ 
one has $e^{-tA} = e^{-t A^\circ} \oplus 0$.
We call $(e^{-tA})_{t > 0}$ the {\bf semigroup generated by $-A$}.

\vertspace

For the remaining part of this paper we fix 
two sectorial forms $\gota$ and $\gotb$ in a Hilbert space $\ch$ 
with $D(\gota) = D(\gotb)$ and vertex~$0$.
It is convenient to shift the index of the $\gota_n$ in Theorem~\ref{tabs103}
by one unit.
For all $n \in \Ni$ define 
\[
\gota_n = \gota + (n-1) \, \gotb
\]
and let $A_n$ be the $m$-sectorial graph associated with $\gota_n$.
Choose $\gotb_1 = \gota$ and $\gotb_n = \gotb$ for all $n \in \Ni$
with $n \geq 2$.
Then it follows from Theorem~\ref{tabsorp301}
that there exists an 
$m$-sectorial graph such that $\lim A_n = A_\infty$ 
in the strong resolvent sense.
Set $A = A_1$, the $m$-sectorial graph associated with $\gota$.

The next theorem is a graph-version of Theorem~\ref{tabs103}\ref{tabs103-1.5}.

\begin{thm} \label{tabs401}
Adopt the above assumptions and notation.
Suppose $\gota$ is closable and there exist $c_1,c_2 > 0$ such that 
\begin{equation}
|\gotb(u)| \leq c_1 \RRe \gota(u) + c_2 \, \|u\|_H^2
\label{etabs401;3}
\end{equation}
for all $u \in D(\gota)$. 
Then there exists an orthogonal projection $P$ in $H$ such that
\[
e^{-t A_\infty}
= \lim_{n \to \infty} \Big( e^{- \frac{t}{n} \, A} \, P \Big)^n
\]
strongly for all $t > 0$.
If in addition $\gota$ is closed, then 
$P$ is the orthogonal projection of $H$ onto 
$ \{ u \in D(\gota) : \gotb(u) = 0 \} \overline{\raisebox{7pt}{$\;\;$}}$, 
where the closure is in $H$.
\end{thm}
\proof\
First assume that $\gota$ is closed.
Then it follows from (\ref{etabs401;3}) that $\gota_n$ is closed for
all $n \in \Ni$.
We use the notation as in Theorem~\ref{tabsorp202}.
Then 
\[
D(\gota_\infty)
= \{ u \in D(\gota) : \gotb(u) = 0 \} 
\]
and $\gota_\infty = \gota|_{D(\gota_\infty) \times D(\gota_\infty)}$.
So by Theorem~\ref{tabsorp202} the graph $A_\infty$ is 
associated with the form~$\gota_\infty$.
Set $W = D(\gota_\infty)$.
Define $\gotq \colon \overline W \times \overline W \to \Ci$ by 
$\gotq(u,v) = 0$ for all $u,v \in \overline W$.
Then $\gotq$ is a closed sectorial form and 
$\gota_\infty = \gota + \gotq$.
Let $Q$ be the graph associated with $\gotq$.
Then $Q = \overline W \times W^\perp$
and $e^{-t Q} = P$ for all $t > 0$.
So by the Trotter--Kato formula, \cite{Kat6} Theorem on page 194, one deduces that 
\[
e^{-t A_\infty}
= \lim_{n \to \infty} \Big( e^{- \frac{t}{n} \, A} \, P \Big)^n
\]
strongly for all $t > 0$.
This proves the theorem if $\gota$ is closed.

If $\gota$ is closable, then apply the above to the closure
of $\gota$ and the extension of $\gotb$ to $D(\overline{\gota})$.\hspace*{5mm}\hfill$\Box$

\vertspace

For the proof of Theorem~\ref{tabs103}\ref{tabs103-2} and the
next example we need a lemma, which allows one to collapse the big direct sum
and the space $W_\infty$ to a subspace of $V_1$ in the proof of Theorem~\ref{tabsorp301}.

\begin{lemma} \label{labs402.5}
Adopt the  assumptions and notation as in the beginning of this section.
Suppose that there exist $c_1,c_2 > 0$ such that 
\[
|\gotb(u)| \leq c_1 \RRe \gota(u) + c_2 \, \|u\|_H^2
\]
for all $u \in D(\gota)$. 
Let $Z$ be the completion of the space $D(\gota)$ with the 
norm $u \mapsto (\RRe \gota(u) + \|u\|_H^2)^{1/2}$.
Let $q \colon D(\gota) \to Z$ be the natural map, 
let $\tilde j \colon Z \to H$ be the extension of the inclusion 
map from $D(\gota)$ into $H$ and let $\tilde \gota,\tilde \gotb \colon Z \times Z \to \Ci$
be the continuous sesquilinear forms such that 
$\tilde \gota(q(u),q(v)) = \gota(u,v)$ and 
$\tilde \gotb(q(u),q(v)) = \gotb(u,v)$ for all $u,v \in D(\gota)$.
Let 
\[
Z_\infty = \{ u \in Z : \tilde \gotb(u) = 0 \}
.  \]
Then $A_\infty$ is the graph associated with 
$(\tilde \gota|_{Z_\infty \times Z_\infty},\tilde j|_{Z_\infty})$.
\end{lemma}
\proof\
For simplicity of presentation we shall consider $Z$ Êas the 
completion of the space $D(\gota)$ with the equivalent norm
$u \mapsto (\RRe \gota(u) + \frac{1}{2} \, \|u\|_H^2)^{1/2}$.
We use the notation as in the proof of Theorem~\ref{tabsorp301}.
Then $Z = V_1$, $q = q_1$, $\tilde j = \Phi_{10}$ and $\tilde \gota = \tilde \gota_1$.
Let $n \in \Ni$.
Then 
\[
\|u\|_n^2 = 2^{-n} \, \|u\|_H^2 + 2^{-(n-1)} \RRe \gota(u) + \sum_{k=2}^n 2^{-(n-k)} \RRe \gotb(u)
\]
for all $u \in D(\gota)$.
Hence 
\[
2^{-(n-1)} \, \|u\|_1^2
\leq \|u\|_n^2
\leq 2 (c_1 + 2 c_2 + 1) \, \|u\|_1^2
\]
for all $u \in D(\gota)$.
Therefore if $1 \leq m \leq n$, then $\Phi_{nm} \colon V_n \to V_m$ is 
bijective.

Let $k,n \in \Ni$ with $2 \leq k \leq n$.
Then $\tilde \gotb_{nk}(q_n(u)) = \gotb(u) = \tilde \gotb( \Phi_{n1}(q_n(u)))$
for all $u \in D(\gota)$.
Hence $\tilde \gotb_{nk}(u) = \tilde \gotb(\Phi_{n1}(u))$ for all $u \in V_n$.
Similarly, $\tilde \gotb_{n1}(u) = \tilde \gota(\Phi_{n1}(u))$ for all 
$n \in \Ni$ and $u \in V_n$.

Let $u \in W_\infty$.
Let $n \in \{ 2,3,\ldots \} $.
Then (\ref{etabsorp301;40}) implies that 
\[
\|\pi_n(u)\|_{V_n}^2
\geq \RRe \tilde \gotb_{nn}(\pi_n(u))
= \RRe \tilde \gotb(\Phi_{n1}(\pi_n(u)))
= \RRe \tilde \gotb(\pi_1(u))
.  \]
Since $\sum_{n=2}^\infty \|\pi_n(u)\|_{V_n}^2 < \infty$ it follows that 
$\tilde \gotb(\pi_1(u)) = 0$ and $\pi_1(u) \in Z_\infty$.
So $\pi_1(W_\infty) \subset Z_\infty$.
We next show that actually an equality is valid.

Let $u \in Z_\infty$.
Then for all $n \in \{ 2,3,\ldots \} $ one has
\begin{equation}
\tilde \gota_n(\Phi_{n1}^{-1}(u))
= \sum_{k=1}^n \tilde \gotb_{nk}(\Phi_{n1}^{-1}(u))
= \tilde \gota(u) + \sum_{k=2}^n \tilde \gotb(u)
= \tilde \gota(u)
.  
\label{elabs402.5;1}
\end{equation}
Hence $(\Phi_{10}(u), u, \Phi_{21}^{-1}(u), \Phi_{31}^{-1}(u), \ldots) \in W_\infty$
by Lemma~\ref{labsorp302}.
Define $\Psi \colon Z_\infty \to W_\infty$ by 
\[
\Psi(u) = (\Phi_{10}(u), u, \Phi_{21}^{-1}(u), \Phi_{31}^{-1}(u), \ldots)
.  \]
Then $\Psi$ is bijective and $\Psi^{-1} = \pi_1|_{Z_\infty}$.

Finally, let $u \in Z_\infty$.
Then $j(\Psi(u)) = \Phi_{10}(u) = \tilde j(u)$
and 
\[
\hat \gota_\infty(u)
= \lim_{n \to \infty} \tilde \gota_n(\pi_n(\Psi(u)))
= \lim_{n \to \infty} \tilde \gota(\Phi_{n1}^{-1}(u))
= \tilde \gota(u)
\]
by (\ref{elabs402.5;1}).
So the graph $A_\infty$ associated (by definition) with $(\hat \gota_\infty, j|_{W_\infty})$
is equal to the graph associated with 
$(\tilde \gota|_{Z_\infty \times Z_\infty},\tilde j|_{Z_\infty})$.\hfill$\Box$

\vertspace

In the next theorem we  assume that (an extension of) the form $\gotb$ 
is associated with a bounded operator.

\begin{thm} \label{tabs402}
Adopt the  assumptions and notation as in the beginning of this section.
Suppose there exists a $($bounded$)$ $B \in \cl(H)$ such that 
\[
\gotb(u,v) = (B u,v)_H
\]
for all $u,v \in D(\gota)$.
Let $P$ be the orthogonal projection from $H$ onto 
$\Big( (B + B^*)(H) \Big)^\perp$.
Then 
\[
e^{-t A_\infty}
= \lim_{n \to \infty} \Big( e^{- \frac{t}{n} \, A} \, P \Big)^n
\]
strongly for all $t > 0$.
\end{thm}
\proof\
We use Lemma~\ref{labs402.5} and its notation.
One has $\tilde \gotb(u,v) = (B (\tilde j(u)),\tilde j(v))_H$ for all $u,v \in Z$.
The operator $B + B^*$ is a positive self-adjoint operator.
Set 
\[
H_1 = \Big( (B + B^*)(H) \Big)^\perp
.  \]
Let $u \in Z$. 
Then $\tilde b(u) = 0$ if and only if $(B + B^*) \tilde j(u) = 0$.
So $\tilde b(u) = 0$ if and only if $\tilde j(u) \in H_1$.
Then $Z_\infty = \{ u \in Z : \tilde j(u) \in H_1 \} $.
Therefore $\ker \tilde j \subset Z_\infty$ and $\ker (\tilde j|_{Z_\infty}) = \ker \tilde j$.
For brevity write
\[
\gott = \tilde \gota|_{Z_\infty \times Z_\infty}
.  \]
Then $A_\infty$ is the graph associated with $(\gott,\tilde j|_{Z_\infty})$.
Define 
\[
V(\gott)
= \{ u \in Z_\infty : \gott(u,v) = 0 \mbox{ for all } v \in \ker (\tilde j|_{Z_\infty}) \} 
\]
(see \cite{AE2} page~36).
Then 
\begin{eqnarray*}
V(\gott)
& = & \{ u \in Z : \tilde j(u) \in H_1 \mbox{ and } \tilde \gota(u,v) = 0 \mbox{ for all } v \in \ker \tilde j \}  \\
& = & V(\tilde \gota) \cap \{ u \in V : \tilde j(u) \in H_1 \} 
,
\end{eqnarray*}
where 
$V(\tilde \gota) 
= \{ u \in Z : \tilde \gota(u,v) = 0 \mbox{ for all } v \in \ker \tilde j \} $.
In particular, $V(\gott) \subset V(\tilde \gota)$
and 
\[
\tilde j(V(\gott)) = \{ \tilde j(u) : u \in V(\tilde \gota) \mbox{ and } \tilde j(u) \in H_1 \}
= \tilde j(V(\tilde \gota)) \cap H_1
.  \]
Define $\tilde \gota_c \colon \tilde j(Z) \times \tilde j(Z) \to \Ci$ by 
$\tilde \gota_c(\tilde j(u),\tilde j(v)) = \tilde \gota(u,v)$ for all $u,v \in V(\tilde \gota)$.
Then $A$ is associated with the closed sectorial form $\tilde \gota_c$ by \cite{AE2} Theorem~2.5(ii).
Similarly, define the form $\gott_c$ with form domain
$D(\gott_c) = \tilde j(Z_\infty)$ by 
$\gott_c(\tilde j(u),\tilde j(v)) = \gott(u,v)$ for all 
$u,v \in V(\gott)$. 
Then $A_\infty$ is the graph associated with $\gott_c$.
Finally, define $\gotq \colon H_1 \times H_1 \to \Ci$ by 
$\gotq(u,v) = 0$ for all $u,v \in H_1$.
Then $\gotq$ is a closed sectorial form and 
$\gott_c = \tilde \gota_c + \gotq$.
Let $Q$ be the graph associated with $\gotq$.
Then $Q = \overline Z_\infty \times Z_\infty^\perp$
and $e^{-t Q} = P$ for all $t > 0$.
Using again the Trotter--Kato formula one deduces that 
\[
e^{-t A_\infty}
= \lim_{n \to \infty} \Big( e^{- \frac{t}{n} \, A} \, P \Big)^n
\]
strongly for all $t > 0$.\hfill$\Box$

\vertspace

The next example show that the boundedness of the operator $B$ in 
Theorem~\ref{tabs402} cannot be replaced by $\gota$-form boundedness of 
the form $\gotb$.

\begin{exam} \label{xabsorp102}
Let $H = L_2(\Ri)$.
Fix $\varphi \in C_c^\infty(\Ri)$ with $\|\varphi\|_H = 1$.
Define $\gota,\gotb \colon W^{1,2}(\Ri) \times W^{1,2}(\Ri) \to \Ci$ by
\begin{eqnarray*}
\gota(u,v) & = & u(0) \, \overline{v(0)} \mbox{ and} \\
\gotb(u,v) & = & (u,v)_H + u(0) \, \overline{v(0)}
   - \tfrac{1}{2} \Big( (u,\varphi)_H + u(0) \Big) \Big( \overline{(v,\varphi)_H} + \overline{v(0)} \Big)
.
\end{eqnarray*}
Then $\gota$ and $\gotb$ are positive symmetric sesquilinear forms.
Moreover, $\gotb(u) \leq 2 \, \gota(u) + 2 \|u\|_H^2$ for all $u \in W^{1,2}(\Ri)$, 
so the form $\gotb$ is $\gota$-bounded.

We follow the construction of the operator $A_\infty$ in Lemma~\ref{labs402.5}.
Define $|||\cdot||| \colon D(\gota) \to [0,\infty)$ by 
$|||u|||^2 = \gota(u) + \|u\|_H^2$.
Define $Z = H \times \Ci$ with the usual inner product and 
define $q \colon D(\gota) \to Z$ by $q(u) = (u,u(0))$.
Then $Z$ is the completion of $(D(\gota),|||\cdot|||)$.
Define $\tilde j \colon Z \to H$ by $\tilde j(u,\lambda) = u$.
Then $\tilde j$ is the continuous extension of the inclusion of $D(\gota)$ into $H$.
The continuous extension $\tilde \gota \colon Z \times Z \to \Ci$ of $\gota$ is 
given by 
\[
\tilde \gota((u_1,\lambda_1),(u_2,\lambda_2))
= \lambda_1 \, \overline{\lambda_2}
.  \]
Let $P_0$ be the orthogonal projection of $Z$ onto $\spann(\varphi,1)$.
Then 
\[
P_0(u,\lambda) = \tfrac{1}{2} ((u,\varphi)_H + \lambda)(\varphi,1)
.  \]
Note that $\gotb(u,v) = ((I - P_0)( q(u) ), q(v))_Z$
for all $u,v \in D(\gota)$.
Hence the continuous extension $\tilde \gotb \colon Z \times Z \to \Ci$ of $\gotb$ is 
given by 
\[
\tilde \gotb((u_1,\lambda_1),(u_2,\lambda_2))
= ((I - P_0)(u_1,\lambda_1),(u_2,\lambda_2))_Z
.  \]
Then 
\[
Z_\infty = \{ v \in Z : \tilde \gotb(v) = 0 \} 
= \spann(\varphi,1)
.  \]
We determine the graph $A_\infty$ associated with 
$(\tilde \gota|_{Z_\infty \times Z_\infty}, \tilde j|_{Z_\infty})$.
Let $(x,f) \in A_\infty$.
Then there exists a $u \in Z_\infty$ such that $x = \tilde j(u)$ and 
$\tilde \gota(u,v) = (f,\tilde j(v))_H$ for all $v \in Z_\infty$.
Let $c \in \Ci$ be such that $u = c \, (\varphi,1)$.
Then $x = c \, \varphi$ and 
\[
(f,\varphi)_H
= \tilde \gota(u,(\varphi,1))
= c \, \tilde \gota((\varphi,1), (\varphi,1))
= c
.  \]
So $f \in c \, \varphi + \varphi^\perp$.
Therefore 
$A_\infty \subset \{ (c \, \varphi, c \, \varphi + y) : c \in \Ci \mbox{ and } y \in \varphi^\perp \} $.
The converse inclusion is easy, so 
\[
A_\infty 
= \{ (c \, \varphi, c \, \varphi + y) : c \in \Ci \mbox{ and } y \in \varphi^\perp \} 
.  \]
Then $e^{-t A_\infty} = e^{-t} \, P_1$ for all $t > 0$, where 
$P_1$ is the orthogonal projection of $H$ onto $\spann \varphi$.

Finally, suppose there exists a projection $P$ such that 
\[
e^{-t A_\infty}
= \lim_{n \to \infty} \Big( e^{- \frac{t}{n} \, A} \, P \Big)^n
\]
strongly for all $t > 0$.
Since $A = 0$ this implies that 
$e^{-t} \, P_1 = \lim_{n \to \infty} P_n = P$
for all $t > 0$.
Therefore $P = P_1 = 0$.
This is a contradiction.
\end{exam}

\subsection*{Acknowledgements}
The second-named author is grateful for a most stimulating stay at the 
University of Oxford.
Part of this work is supported by an
NZ-EU IRSES counterpart fund from Government funding,
administered by the Royal Society of New Zealand.


\begin{thebibliography}{AEKS13}

\bibitem[AE]{AE2}
{\sc Arendt, W. {\rm and} Elst, A. F.~M. ter}, Sectorial forms and degenerate
  differential operators.
\newblock {\em J. Operator Theory} {\bf 67} (2012),  33--72.

\bibitem[AEKS]{AEKS}
{\sc Arendt, W., Elst, A. F.~M. ter, Kennedy, J.~B. {\rm and} Sauter, M.}, The
  Dirichlet-to-Neumann operator via hidden compactness, 2013. Preprint.
\newblock ArXiv:1305.0720.

\bibitem[Br{\'e}]{Bre2}
{\sc Br{\'e}zis, H.}, {\em Op\'erateurs maximaux monotones et semi-groupes de
  contractions dans les espaces de Hilbert}.
\newblock North-Holland Mathematics Studies 5. North-Holland Publishing Co.,
  1973.
\newblock Notas de Mathem\'atica (50).

\bibitem[FW]{FW}
{\sc Fillmore, P.~A. {\rm and} Williams, J.~P.}, On operator ranges.
\newblock {\em Adv. Math.} {\bf 7} (1971),  254--281.

\bibitem[Kat1]{Kat6}
{\sc Kato, T.}, Trotter's product formula for an arbitrary pair of self-adjoint
  contraction semigroups.
\newblock In {\em Topics in functional analysis}, Adv. in Math. Suppl. Stud. 3.
  Academic Press, New York, 1978,  185--195.

\bibitem[Kat2]{Kat1}
\leavevmode\vrule height 2pt depth -1.6pt width 23pt, {\em Perturbation theory
  for linear operators}.
\newblock Second edition, Grund\-lehren der mathematischen Wissenschaften 132.
  Springer-Verlag, Berlin etc., 1980.

\bibitem[Ouh1]{Ouh7}
{\sc Ouhabaz, E.-M.}, Second order elliptic operators with essential spectrum
  $[0,\infty)$ on $L^p$.
\newblock {\em Comm. Partial Differential Equations} {\bf 20} (1995),
  763--773.

\bibitem[Ouh2]{Ouh5}
\leavevmode\vrule height 2pt depth -1.6pt width 23pt, {\em Analysis of heat
  equations on domains}, vol.\ 31 of London Mathematical Society Monographs
  Series.
\newblock Princeton University Press, Princeton, NJ, 2005.

\bibitem[Sho]{Sho}
{\sc Showalter, R.~E.}, {\em Monotone operators in Banach space and nonlinear
  partial in differential equations}.
\newblock Mathematical Surveys and Monographs 49. American Mathematical
  Society, Providence, RI, 278.

\bibitem[Sim]{bSim5}
{\sc Simon, B.}, A canonical decomposition for quadratic forms with
  applications to monotone convergence theorems.
\newblock {\em J.\ Funct.\ Anal.} {\bf 28} (1978),  377--385.

\end{thebibliography}
\end{document}